\documentstyle[12pt]{article}
\begin{document}

\begin{titlepage}

August 10th 1999 \hfill
\vskip 1.5cm
\centerline{\large \bf
  Universal R-matrix formalism for the spin Calogero-Moser }
\centerline{\large \bf system and its difference counterpart}
\vskip 1.0cm
\centerline{D. Talalaev
\footnote{E-mail:talalaev@css-mps.ru} 
}
\centerline{High Geometry and Topology Department, Faculty of Mathematics
and Mechanics,}
\centerline{Moscow State University
\footnote{M.V. Lomonosov Moscow State University,
Vorobjevy Gory, Moscow, 119899, Russia.}}
\vskip 3.0cm
{\bf Abstract.} The expression of the quantum Ruijsenaars-Schneider
Hamiltonian is obtained in the framework of the dynamical $R$-matrix
formalism. This generalizes to the case of $U_q(sl_n)$ the result obtained
in \cite{BBB} for $U_q(sl_2)$ which is the higher difference Lame operator.
The same result was
obtained in \cite{FelVar} using the ``fusion'' procedure. The distinctive
starting point of the present paper is the universal dynamical $R$-matrix.

\end{titlepage}

\section{Introduction}
The connection between the finite gap integrability and the R-matrix
formalism comes from the famous work \cite{TaFa} and the claim therein that
the Lax representation 
$$\dot L =[L,M]$$
is equivalent to the existence of a classical r-matrix satisfying the
following fundamental property:
$$ \{L_1,L_2\}=[r_{12},L_1]-[r_{21},L_2];$$
$r$ lies in the tensor product $Mat(n)^{\otimes 2}$, where $L_1=L \otimes
id$.
This sets the path to quantization. For example, in \cite{ABB}, it was
shown that in the case of the spin Calogero-Moser model one can treat the
dynamical Yang-Baxter equation (Gervais-Neveu-Felder equation \cite{GN}) as
the quantized r-matrix representation for this model. There the GNF
equation in the form of Gervais-Neveu writes:
$$R_{12}(x)R_{13}(xq^{H_2})R_{23}(x)=
R_{23}(xq^{H_1})R_{13}(x)R_{12}(xq^{H_3})$$
where $R(x)$ lies universally in $U_q(sl_n) \otimes U_q(sl_n)$ or in the
tensor product of the representations. 
Also is known the L-variant form of this equation:
$$R_{12}(xq^{-H_3 /2})L_{13}(x)L_{23}(x)=
L_{23}(x)L_{13}(x)R_{12}(xq^{H_3 /2}).$$
The formula above provides a machinery for the construction of the quantum
integrable system, as was done in \cite{ABB}. The main claim of this
procedure is that the operators $I_n$ on some space of functions with
values in some fixed representation of $U_q(sl_n)$ constructed by the
formula
$$I_n=Tr_{1...n}[L_1(x)...L_n(x) \hat R_{12}(xq^{h^{(3,n)}})...\hat
R_{n-1,n}(x)] $$
where $\hat R=PR$, $P$ is the permutation operator in the tensor product
and $h^{(k,l)}=\sum_{i=k}^l h^{(i)}$,  
leave the subspace of zero weight in the representation invariant. The
restrictions of the operators $I_n$ to this subspace form a set of
commuting operators.

The crucial ideas to obtain the integrable system are:

1. take the solution of the standard Yang-Baxter equation as the universal
Drinfeld R-matrix \cite{A} for $U_q(sl_n)$;

2. use the quasi-Hopf shift $F_{12}(x)$ \cite{J} to obtain the quasi-Hopf
structure on this algebra and the solution of the GNF equation in the form 
$F_{21}^{-1}(x)R_{12}F_{12}(x)$;

3. take as the L-operator the expression 
\begin{equation}
L=q^{-(H_1+H_2 /2)p} R q^{H_2 p/2},
\label{eq:L}
\end{equation}
 where $p=x \frac {\partial} {\partial x}$ in the tensor product of the
representatian $\rho \otimes \rho_q$, where $\rho$ is the matrix
representation and $\rho_q$ is some quantum representation;

4. the set of commuting $I_n$ operators form an integrable system.

This scheme was successively investigated in \cite{BBB} for the case
$U_q(sl_2)$. For the spin-$j$ representation on the zero-weight subspace
the element $I_1$ was shown to be the Hamiltonian generalizing, in the case
of spin-j, the one-particle Ruijsenaars-Schneider system. The quantum
problem was solved  using the existence of a spin-shift operator. The
constructed quantum Hamiltonians were shown to be, in the limit $q \mapsto
1$, those of the Calogero-Moser system.

The main goal of the present paper is to generalize this construction for
the case $U_q(sl_n)$. This result can be treated as en exercise but the
evident 
advantage of this approach is that it allows one to investigate the full
representation, full algebra and not only the zero-weight space. 
It is likely to clarify the original Poisson structure of the spin RS
system.

In the second part, I explain the method to construct the Drinfeld
R-matrix,
emphasizing the generalization to $U_q(sl_n)$.

In the third part I explicitly calculate the dynamical $R$-matrix by using
the linear equation satisfied by the quasi-Hopf shift.

In the fourth and last part I construct the Hamiltonian of the system.

\section{Universal R-matrix.}
The construction of the universal R-matrix arises from the notion of
Drinfeld double for quasitriangular Hopf algebras. It is the case for
$U_q(sl_n)$, described by the set of generators
$f_\alpha,e_\alpha,h_\alpha$; for primitive roots $\alpha_i$ we shall write
them for simplicity $f_i,e_i,h_i$, $i=1,...,n-1$. They satisfy the
following relations:
$$[h_i,h_j]=\delta_{ij},\qquad [e_i,f_j]=\delta_{ij}\frac
{q^{h_i}-q^{-h_i}}  {q-q^{-1}}$$
$$[h_i,e_j]=(\delta_{ij}2-\delta_{i,j\pm 1})e_j, \qquad
  [h_i,f_j]=(-\delta_{ij}2+\delta_{i,j\pm 1})f_j.$$
The Serre's relations in this case read:
$$[e_i,e_j]=0 \qquad if |i-j|>1 \qquad [f_i,f_j]=0 \qquad if |i-j|>1$$
$$e_ie_ie_{i+1}+e_{i+1}e_ie_i=(q+q^{-1})e_ie_{i+1}e_i,$$
and the same for the $f_i$'s.

For a fixed order of the positive roots one introduces canonical generators
by requiring that, for $\alpha, \beta,\alpha+\beta$ - positive roots such
that $\alpha<\alpha+\beta<\beta$, we have 
$$e_{\alpha+\beta}=e_\alpha e_\beta - q e_\beta e_\alpha \qquad
f_{\alpha+\beta}=f_\beta f_\alpha - q^{-1} f_\alpha f_\beta.$$

We associate the elements $e_{i,i+1},f_{i,i+1}$ to the primitive generators
$e_i,f_i$ respectively, and extend this by induction as 
$$e_{ij}=e_{ik} e_{kj} - q e_{kj} e_{ik}, \qquad
f_{ij}=f_{kj} f_{ik} - q^{-1} f_{ik} f_{kj} \qquad for\ i<k<j.$$
Due to the Serre's relations it does not depend on the choice of k. The
order here is lexicographical, like $(12),(13),...(1,n),(2,3),...(n-1,n)$.

A Poincare-Birkhoff-Witt basis of $U_q(N_+)$ (resp.$U_q(N_-))$ is formed
here by:
\begin{equation}
e^p=\prod_{\alpha \in \Phi^+}^>(e_{ij})^{p_{ij}} \qquad
(\mbox{resp.} f^p=\prod_{\alpha \in \Phi^+}^>(f_{ij})^{p_{ij}}).
\label{eq:PBW}
\end{equation}

The R-matrix can be expressed as:
$$R=K \hat R, \qquad K=\prod_{j=1}^{n-1} q^{h_i \otimes h^i}, \qquad
\hat R = \prod_{\Phi^+,>} \hat R_\alpha,$$
where $\hat R_\alpha=\exp_q((q-q^{-1})e_\alpha \otimes f_\alpha)$. The
q-exponential is defined as follows:
$$\exp_q(z)=\sum_{n=0}^\infty \frac 1 {[n]!} z^n,$$
where $[n]=(q^n-q^{-n})/(q-q^{-1}).$

We will be interested in calculating the R-matrix when the first or the
second space is the space of matrix representation of $U_q(sl_n)$; then, in
the series for the q-exponential, there is no terms with degrees of
generators but 0 and 1, i.e. the formula for the R-matrix reads:
$$\hat R=\prod_{ij,>}(1+(q-q^{-1})e_{ij} \otimes f_{ij}).$$
 We will separately treat two cases:

1. matrix representation on the first space: $e_{ij}$ is just the matrix
element $\delta_{ij}$ for $i<j$. Looking at the order in the product and
noting that $\delta_{ij}\delta_{km}=0$ if $(ij)>(km)$ we have
$$\hat R=1+(q-q^{-1})\sum_{i<j}e_{ij}\otimes f_{ij};$$

2. matrix representation on the second space: $f_{ij}$ is then the matrix
element $\delta_{ji}$ for $i<j$ and we have for the R-matrix:
$$\hat R=1+\sum_{i<j} (\sum
e_{i_{m-1}i_m}...e_{i_1i_2}(q-q^{-1})^{m-1})\otimes f_{ij}.$$
The inner summation is taken on all partitions $i=i_1<i_2<...<i_m=j$.

\vskip 1.5cm
\section{Quasi-Hopf shift.}
The main representation of $U_q(sl_n)$ on the quantum space is a highest
weight representation (nl,0,...,0). In the latter the weight spaces have
the fundamental property of being one-dimensional. For example on the
zero-weight space all the elements $f_{i_{m-1}i_m}...f_{i_1i_2}$ with the
same $i_1,i_m$ are proportional. Looking for the representations of
$U_q(sl_2)$ subalgebras corresponding to all positive roots we find
$$f_{i_{m-1}i_m}...f_{i_1i_2}v_0
=([l+1]q^{-l})^{(m-1-i_m+i_1)}f_{i_1i_m}^*v_0$$
where $f_{ij}^*=f_{j-1,j}f_{j-2,j-1}...f_{i,i+1}$ and $v_0$ is the
zero-weight vector. Introducing the element 
$e_{ij}^*=e_{j-1,j}e_{j-2,j-1}...e_{i,i+1}$ 
we have on the space generated by
$f_{ij}^* v_0$
$$ e_{i_{m-1}i_m}...e_{i_1i_2}v_0
=([l]q^{l})^{(m-1-i_m+i_1)}e_{i_1i_m}^*.$$
It is worth saying that those relations are preserved through respectively
right and left multiplication by a positive root generator $e_{i_mj}$ and
$f_{i_mj}$. It will be important for further calculations.

This structure of representation leads for example to the property that, in
the case of matrix representation on the second space, the R-matrix reads:
\begin{equation}
\hat R_{21}=1+\sum_{i<j} f_{ij}\otimes e_{ij}^*(q^l [l]^{-1})^{j-i-1}. 
\label{eq:R}
\end{equation}
The last fact needed from the representation is
\begin{equation}
e_{ij}^*f_{ij}^*v_0=([l][l+1])^{j-i}v_0.
\label{eq:ef}
\end{equation}
The construction of the quasi-Hopf $F_{12}(x)$ shift is based on the main
proposition of \cite{A}, which claims that it satisfies the linear
equation:
$$F_{12}(x)B_2(x)=\hat R_{12}^{-1}B_2(x)F_{12}(x)$$
or its equivalent form:
\begin{equation}
\hat R_{12}F_{12}(x)= B_2(x)F_{12}(x) B_2(x)^{-1}
\label{eq:RF}
\end{equation}
where $B(x)=q^{\sum_{j=1}^{n-1}(h_jh^j-x_jh^j)}$. The element $F_{12}(x)$
must be of the form:
$$F(x)=\sum_{p,r}e^p \otimes f^r \phi_{p,r}(x)$$
where the summation is done over the Poincare-Birkhoff-Witt basis elements 
(\ref{eq:PBW}) with the same weight. Introducing the variables $y_j$ such
that $x_i=y_{i+1}-y_i$ we have 
\hbox{$Bf_{ij}B^{-1}=f_{ij}q^{y_j-y_i-2h_{ij}+2}$} where
$h_{ij}=\sum_{k=i}^{<j}h_k$. Then, taking $F$ of the form:
$$F_{12}(x)=1+\sum_{i<j}e_{ij} \otimes f_{ij}^* \Phi_{ij}(x)$$
one finds
$$B_2 F_{12}B_2^{-1}=1+\sum_{i<j}e_{ij} \otimes f_{ij}^*
\Phi_{ij}(x)q^{y_j-y_i+2},$$
where we have restricted the right-hand-side of the equation above to the
zero-weight subspace.
Treating alike its left-hand-side one obtains:

\vbox{
$$\hat R_{12}F_{12}(x)=1+\sum_{i<j}e_{ij}\otimes(f_{ij}^*\Phi_{ij}+
f_{ij}(q-q^{-1})$$
$$+ (q^2-1)\sum_{k,i<k<j}f_{j-1,j}...f_{k+1,k+2}
(f_{k,k+1}f_{ik}-f_{i,k+1})\Phi_{kj}).$$
}

Reducing the proportional terms we have:
$$\hat R_{12}F_{12}(x)=1+\sum_{i<j}e_{ij}\otimes f_{ij}^*(\Phi_{ij}
+q_l^{i-j+1}(q-q^{-1})+(q^2-1)\sum_{k,i<k<j}q_l^{i-k}(q_l-1)\Phi_{kj}),$$
where $q_l=[l+1]q^{-l}$.
Comparing the two sides of (\ref{eq:RF}) that we just evaluated yields an
equation to be satisfied by the functions $\Phi_{ij}$:
$$\Phi_{ij}(q^{y_j-y_i+2}-1)=q_l^{i-j-1}(q-q^{-1})+
\sum_{k,i<k<j}(q^2-1)(q_l-1)q_l^{i-k}\Phi_{kj}.$$
It leads to the relations:
\begin{equation}
\Phi_{ij}=\Phi_{i+1,j}q_l^{-1}\frac {q^{y_j-y_{i+1}+2}-q^{-2l}} 
{q^{y_j-y_i+2}-1}, \qquad \Phi_{j-1,j}=\frac {q-q^-1}
{q^{y_j-y_{j-1}+2}-1}.
\label{eq:Phi}
\end{equation} 
To evaluate the element $F_{21}^{-1}$ we rewrite (\ref{eq:RF}) as follows:
\begin{equation}
F_{21}^{-1}=B_1 F_{21}^{-1}B_1 \hat R_{21}.
\label{eq:F}
\end{equation}
Taking then $F_{21}^{-1}$ in the form:
$$ F_{21}^{-1}=1+\sum_{i<j}\Psi_{ij}f_{ij}\otimes e_{ij}^* $$
we get
$$ B_1 F_{21}^{-1}B_1^{-1}=
1+\sum_{i<j}\Psi_{ij}q^{y_j-y_i}f_{ij}\otimes e_{ij}^* ,$$
and, plugging the expression (\ref{eq:R}) for $R_{21}$ in (\ref{eq:F})
gives

\vbox{
$$ B_1 F_{21}^{-1}B_1^{-1}\hat R_{21}=1+\sum_{i<j}f_{ij}\otimes e_{ij}^*
(\Psi_{ij}q^{y_j-y_i}+(q-q^{-1}) (q^l [l]^{-1})^{j-i-1} $$
$$+\sum_{k,i<k<j}\Psi_{kj}q^{y_j-y_k}(q-q^{-1}) (q^l [l]^{-1})^{k-i-1}) .$$
}

Identifying the two sides of (\ref{eq:F}) gives a system of equations for
the functions $\Psi_{ij}$:
$$\Psi_{ij}(1-q^{y_j-y_i})=(q-q^{-1})(q^l [l]^{-1})^{j-i-1} +
\sum_{k,i<k<j}\Psi_{kj}q^{y_j-y_k}(q-q^{-1}) (q^l [l]^{-1})^{k-i-1}).$$
Solving this system we get:
\begin{equation} 
\Psi_{ij}=q^{-l}[l]^{-1} \frac {q^{y_j-y_{i+1}}-q^{2l}} {q^{y_j-y_i}-1}
\Psi_{i+1,j}
\qquad \Psi_{j-1,j}=\frac {q-q^{-1}} {1-q^{y_j-y_{j-1}}}.
\label{eq:Psi}
\end{equation}

\section{ The Hamiltonian.}
In this last step all calculations are performed in the matrix
representation on the first space and in the zero weight subspace of the
(nl,0,...,0)-highest weight representation on the quantum space. We are
going to calculate the product:
$$F_{21}^{-1}K \hat R_{12}F_{12}$$
where 
$$K \hat R_{12} F_{12}(x)=1+\sum_{i<j}e_{ij} \otimes f_{ij}^*
\Phi_{ij}(y)q^{y_j-y_i+1}. $$
Because we will take traces, only the diagonal elements of $R_{12}(x)$
matter. They can be expressed as:
$$R_{12}(y)_{jj}=1+\sum_{i<j}e_{ij}^*f_{ij}^*
\Phi_{ij}(y)\Psi_{ij}(y)q^{y_j-y_i+1},$$
which can be simplified by using equation (\ref{eq:ef}) into:
$$R_{12}(y)_{jj}=1+\sum_{i<j}\Phi_{ij}(y)
\Psi_{ij}(y)q^{y_j-y_i+1}([l][l+1])^{j-i}.$$
The reduced form of this sum is
$$R_{12}(y)_{jj}^k=
1+\sum_{i=k}^{i<j}\Phi_{ij}(y)\Psi_{ij}(y)q^{y_j-y_i+1}([l][l+1])^{j-i}.$$
{\bf Lemma 1.} For $\tilde R_j^k= R_{12}(y)_{jj}^k$ the following property
is fulfilled:
$$\tilde R_j^k=\prod_{i=k}^{i<j}
(1-[l][l+1]\xi(y_j-y_i))$$
where $$\xi(y)=\frac {(q-q^{-1})^2 q^{y+1}} {(q^y-1)(q^{y+2}-1)}.$$

{\it Proof.} For $k=j-1$ it is true. We then proceed by induction.
Supposing the claim valid for $k<m$ leads to 
$$\tilde R_j^{k}=\tilde R_j^{k+1}(1-[l][l+1]\xi(y_j-y_k)),$$
and it is then sufficient to prove that 
\begin{equation}
\Phi_{kj}\Psi_{kj}q^{y_j-y_k+1}([l][l+1])^{j-k}=
-[l][l+1]\xi(y_j-y_k) \prod_{i=k+1}^{i<j}(1-[l][l+1]\xi(y_j-y_i)).
\label{eq:ind2}
\end{equation}
Carrying out the induction for (\ref{eq:ind2}) requires proving
$$ \frac {\Phi_{k-1,j}}{\Phi_{kj}} \frac {\Psi_{k-1,j}} {\Psi_{kj}}
q^{y_k-y_{k-1}}=\frac {\xi(y_j-y_{k-1})} {\xi(y_j-y_k)}
(1-[l][l+1]\xi(y_j-y_k)),$$
which is straightforwardly achieved by plugging in it the expressions for
the $\Phi$'s and $\Psi$'s obtained in (\ref{eq:Phi}) and (\ref{eq:Psi}).
This proves  lemma 1.
\vskip .5cm
{\bf Lemma 2.} The Hamiltonian reads:
$$H=I_1=\sum_{i=1}^n e^{p_i} \prod_{k<i}(1-[l][l+1]\xi(y_i-y_k))$$
where $p_i=2 \frac \partial {\partial_{y_i}}$

{\it Proof.} The main thing to prove is that 
$q^{pH}=q^{\sum \tilde p_i h^i}$  of (\ref{eq:L}) in the matrix
representation is just the diagonal
matrix with elements $e^{p_i}$, $p_i$ as above. The variables used in
\cite{BBB} are $\tilde {x_i}=q^{-\frac {y_{i+1}-y_i} 2}$. There, the
authors introduce the operator $\tilde p_i=\tilde x_i \frac {\partial}
{\partial_{\tilde x_i}}.$  The link between the $p_i$'s of this paper and
the $\tilde p_i$'s is $q^{\tilde p_i}=e^{p_{i+1}-p_i}$. In this
representation the elements $h^{i-1}-h^i$, $h^0=h^n=0$, are just the
diagonal $\delta_{ii}$'s. So 
$$q^{\sum \tilde p_i h^i}=e^{\sum h^i (p_{i+1}-p_i)}
=e^{\sum p_i (h^{i-1}-h^i)}=\sum \delta_{ii} e^{p_i}.$$
Introducing the variables $t_i=y_i/2,t=y/2$ we have 
$$\xi(y)=\frac {(q-q^{-1})^2} {(q^t-q^{-t}) (q^{t+1}-q^{-t-1})}.$$
Conjugating the Hamiltonian
$$\hat H_l=f H_l f^{-1}$$
where the function $f(t_1,..,t_n)=\prod_{i<k}g(t_i-t_k)$ and 
$g(t)=((q^{t+1}-q^{-t-1})...(q^{t+l}-q^{-t-l}))^{-1}$ 
one obtains
$$\hat H_l=\sum_i e^{\frac \partial \partial_{t_i}}
\prod_{k \ne i} \frac {q^{t_i-t_k+l}-q^{-t_i+t_k-l}}
{q^{t_i-t_k}-q^{-t_i+t_k}}.$$

The elliptic version of this Hamiltonian was obtained in \cite{FelVar}. The
similar Hamiltonian 
$$H_{\pm}=\sum_i (e^{i\theta_i} \pm e^{-i\theta_i})
\prod_{j \ne i}f(q_{ij}),$$
where
$$f^2(q)=(1-\frac {sin^2(\pi\nu/kN)} {sin^2(\pi q/k)}),$$
was obtained in \cite{Gor} in the framework of hamiltonian reduction for
the group $G$ on the 
space $T^*\hat G$ where $\hat G$ is the suitable affine Lie group. 

\vskip 1.0cm
{\em Acknowledgements}

I would like to thank my scientific advisor I. Krichever for introducing me
to 
the beauty of integrable systems. 

I would like then to thank LPTHE which welcomed me for its support, the
friendly and studious atmosphere,
that I found there, specially O. Babelon for leading me into this problem
and for his precious advice and encouragements.

\end{document}